\documentclass[amsfonts,12pt]{article}
\usepackage{amssymb}

\newcommand{\sll}[2]{{\rm SL}({#1},{#2})}
\newcommand{\gl}[2]{{\rm GL}({#1},{#2})}

\newcommand{\seisei}[1]{\mathop{<}{#1}\mathop{>}}
\newcommand{\oku}{\mathrel{:=}}

\newcommand{\sing}{{\rm Sing}}
\newcommand{\spec}{{\rm Spec}}

\newtheorem{Theorem}{Theorem}[section]

\newtheorem{Lemma}[Theorem]{Lemma}

\newcommand{\qed}{{\unskip\nobreak\hfil\penalty50\quad\null\nobreak\hfil{\bf
Q.E.D.}\parfillskip0pt\finalhyphendemerits0\par\medskip}}

\newcommand{\proof}{\noindent{\bf Proof.} \ }

\newcommand{\rmn}[1]{\romannumeral#1}
\def\RMN#1{\uppercase\expandafter{\romannumeral#1}}

\begin{document}

\title{Gorenstein isolated quotient singularities 
of odd prime dimension are cyclic}
\date{}
\author{Kazuhiko Kurano and Shougo Nishi}
\maketitle

\abstract{In this paper, we shall prove that 
Gorenstein isolated quotient singularities 
of odd prime dimension are cyclic.
In the case where the dimension is bigger than $1$ and is not 
an odd prime number, then there exist
Gorenstein isolated non-cyclic quotient singularities.
}

\vspace{3mm}


\vspace{3mm}

\noindent
{\bf  2000 Math.\ Subject Classification}.  primary 13A50, secondary 14L30

\section{Introduction}\label{Intro}

Let $G$ be a finite subgroup of $\gl{n}{{\Bbb C}}$,
where ${\Bbb C}$ is the field of complex numbers and let $\gl{n}{{\Bbb C}}$
be the set of $n \times n$ invertible matrices with entries in ${\Bbb C}$.
Then, $G$ acts on a polynomial ring $R = {\Bbb C}[X_1,X_2,\ldots, X_n]$
linearly.
Let $R^G$ be the invariant subring, i.e., 
\[
R^G = \{ r \in R \mid g(r) = r \ \ \forall g \in G \} .
\]
It is well-known that $R^G$ is finitely generated over ${\Bbb C}$
(cf.\ Theorem~1.3.1 in \cite{Benson}).

It is possible to classify finite subgroups in $\sll{2}{{\Bbb C}}$
(cf.\ Theorem~2.4.5 in \cite{Watanabe1}).
Here, $\sll{n}{{\Bbb C}}$ is the subgroup of $\gl{n}{{\Bbb C}}$
consisting of all matrices of determinant $1$.
It is well-known that 
the invariant subring of ${\Bbb C}[X_1,X_2]$  under the linear action
of a finite subgroup of $\sll{2}{{\Bbb C}}$
is a hypersurface in ${\Bbb C}^3$ with isolated singularity.

It is also possible to classify finite subgroups in $\sll{3}{{\Bbb C}}$ 
(cf.\ Yau-Yu~\cite{YY}).
Using the classification, it was proven that
Gorenstein isolated quotient singularities 
of dimension three are cyclic (Theorem~A and Theorem~23 in Yau-Yu~\cite{YY}).

The purpose of  this paper is to prove the following theorem:

\begin{Theorem}\label{mainc}
Let $n$ be an odd prime number.
Let $G$ be a finite subgroup of $\gl{n}{{\Bbb C}}$ which  contains
no pseudo-reflection.
Assume that the invariant subring $R^G$ is Gorenstein 
with isolated singularity.
Then, $R^G$ has a cyclic quotient singularity.
\end{Theorem}

For a finite subgroup $G$ of $\gl{n}{{\Bbb C}}$, we set  
\[
\Sigma_i=\{g\in G \mid \mbox{$1$ is an eigenvalue of $g$ 
with multiplicity at least $i$} \}
\]
for $i=0, 1, \ldots, n$.
Each element in $\Sigma_{n-1}\backslash \{e\}$ is called a 
{\em pseudo-reflection}.
Set
\[
H_i=\seisei{\Sigma_{i}} ,
\]
which is the subgroup of $G$ generated by $\Sigma_{i}$.
By definition we have 
\[
\begin{array}{l}
G=\Sigma_0 \supset \Sigma_1 \supset \cdots \supset \Sigma_{n-1} \supset \Sigma_n=\{e\} \ \ \mbox{and}  \\
G=H_0 \supset H_1 \supset \cdots \supset H_{n-1} \supset H_n=\{e\}  .
\end{array}
\]
Here, remark that $\Sigma_n$ is equal to $\{ e\}$,
since any matrix in $G$ is diagonalizable.
These are very important subgroups, because
the ring homomorphism $R^G \rightarrow R^{H_l}$ is \'etale in codimension $s$
if and only if $l \le n - s$.

Suppose $n \geq 2$.
Let $l$ be an integer such that $0\leq l \leq n-2$.
By purity of branch locus (cf.\ Theorem~41.1 in \cite{Nagata}) and 
the Shephard-Todd theorem (cf.\ Theorem~7.2.1 in \cite{Benson}), 
we know that the following two conditions are equivalent:
\begin{description}
\item[$(1)$]
$H_l\supsetneqq H_{l+1}=\cdots =H_{n-1}$, 
\item[$(2)$]
$\sing R^G \neq \o$ and $\dim \sing R^G=l$.
\end{description}
Here $\sing R^G$ is the {\em singular locus} of $R^G$, i.e.,
\[
\sing R^G = \{ P \in \spec R^G \mid 
\mbox{$(R^G)_P$ is not a regular local ring}
\} .
\]
If $\sing A$ is not empty and if the dimension of $\sing A$ is $0$,
we say that $A$ has {\em isolated singularities}.
Thus, the following two conditions are equivalent:
\begin{description}
\item[$(1)$]
$R^G$ has isolated singularities.
\item[$(2)$]
$H_0\supsetneqq H_1=\cdots =H_{n-1}$D
\end{description}
If $\Sigma_{n-1} =\{e\}$, then the above two conditions are equivalent to
the following:
\begin{description}
\item[$(3)$]
$\Sigma_1=\{e\}$, i.e., $1$ is not an eigenvalue of 
any matrix in $G$ except for $e$.
\end{description}

On the other hand, remember the following theorem due to Watanabe~\cite{Watanabe}:

\begin{Theorem}[Watanabe]\label{wata}
Let $G$ be a finite subgroup of $\gl{n}{{\Bbb C}}$ and
suppose that $G$ acts on $R\oku {\Bbb C}[X_1,X_2,\ldots ,X_n]$ linearly.
\begin{itemize}
\item
If $G\subset \sll{n}{K}$, then $R^G$ is a Gorenstein ring.
\item
If $R^G$ is a Gorenstein ring and if $\Sigma_{n-1} =\{e\}$,
then $G\subset \sll{n}{K}$.
\end{itemize}
\end{Theorem}

Since $R^{H_{n-1}}$ is isomorphic to a polynomial ring and 
$G/H_{n-1}$ acts on $R^{H_{n-1}}$ linearly,
the case where $\Sigma_{n-1}=\{e\}$ is very important.

By these arguments, 
if $\Sigma_{n-1}=\{e\}$, we have the following assertions:
\begin{itemize}
\item
$G\subset \sll{n}{K}$ if and only if $R^G$ is Gorenstein.
\item
$R^G$ has isolated singularities
if and only if $1$ is not an eigenvalue of any matrix in $G$ except for $e$.
\end{itemize}

Thus Theorem~\ref{mainc} immediately follows from Lemma~\ref{main} below.

\begin{Lemma}\label{main}
Let $n$ be an odd prime number.
Let $G$ be a finite subgroup of $\sll{n}{K}$, 
where $K$ is a field such that the characteristic of $K$ is $0$ or 
does not divide the order of $G$.
Assume that 
$1$ is not an eigenvalue of any matrix in $G$ 
except for the unit matrix.
Then, $G$ is a cyclic group.
\end{Lemma}

We remark that the pair $(G, \rho)$ of a finite group $G$ and its irreducible 
fixed point free complex representation $\rho$ are classified,
where fixed point free means that $\rho(s)$ does not have $1$ as its eigenvalue for $s \neq e$.
This classification is obtained in Theorem~7.2.18 in \cite{wolf}.
Therefore, Lemma~\ref{main} follows from the classification.

In this paper, we give a very simple and elementary proof to Lemma~\ref{main}.

We shall prove Lemma~\ref{main} in Section~\ref{proof}.
In Section~\ref{rei}, we shall give examples of non-cyclic subgroups in 
the case where $n$ is bigger than $1$ and is not an odd prime integer.

\section{Proof of Lemma~\ref{main}}\label{proof}

We shall prove Lemma~\ref{main} in this section.

We may assume that $K$ is an algebraically closed field.

Remark that each matrix in $G$ is diagonalizable
because the characteristic of $K$ is $0$ or
does not divide the order of $G$.

First we shall prove Lemma~\ref{main} in the case where
$G$ is an abelian group.
Next we shall do in the case where $G$ is a solvable group.
Finally we prove Lemma~\ref{main} 
without any additional assumptions.

\subsection{The case where $G$ is abelian}\label{abel}

In this subsection, we prove Lemma~\ref{main} in the case where
$G$ is an abelian group.

Assume that $G$ is a finite abelian subgroup of $\sll{n}{K}$D

Since the characteristic of $K$ is $0$ or does not divide the order of $G$,
there exists $c\in \gl{n}{K}$ such that $c^{-1}gc$ is a diagonal matrix
for any $g\in G$.
Set $c^{-1}Gc\oku \{ c^{-1}gc | g\in G\}$.
Remember that  $g$ and $c^{-1}gc$ have the same characteristic polynomial.
So, $g$ and $c^{-1}gc$ have the same determinant and the same eigenvalues.
Replacing $G$ with $c^{-1}Gc$, we may assume that 
all matrices in $G$ are diagonal.

We define
$$ \psi : G \longrightarrow K^\times $$
by letting $\psi(g)$ be the $(1,1)$th entry of each diagonal matrix $g$ in $G$.
Then, it is a group homomorphism.
Since $1$ is not an eigenvalue of any matrix in $G$ except for 
the unit matrix,
$\psi$ is injective.

Since any finite subgroup of $K^\times$ is cyclic,
so is $G$.

\subsection{The case where $G$ is solvable}\label{skakai}

In this subsection, we prove Lemma~\ref{main} in the case where
$G$ is a solvable group by induction on ${}^\# G$ (the order of $G$).

Let $G$ be a finite solvable subgroup of $\sll{n}{K}$
satisfying the assumption in Lemma~\ref{main}D
Assume ${}^\# G > 1$.
By induction, any finite solvable subgroup $G'$ of $\sll{n}{K}$
satisfying the assumption in Lemma~\ref{main} is cyclic
if ${}^\# G > {}^\# G'$.
In particular, any proper subgroup of $G$ is cyclic.

Let $H$ be a maximal subgroup of $G$ 
that contains the commutator subgroup of $G$.
We remark that such a subgroup exists since $G$ is solvable.
Then $H$ is a normal subgroup of $G$.
Since $H$ is a proper subgroup of $G$,
$H$ is a cyclic group.
Let $a$ be a generator of $H$, and take $b \in G \setminus H$.
Then, 
\[
\mbox{$H = \seisei{a}$ and $G=\seisei{a,b}$},
\]
where $\seisei{a_1, \ldots, a_t}$ means the subgroup generated by 
$a_1$, \ldots, $a_t$.

Let $s$ be the order of $a$.
Since $H$ is a normal subgroup of $G$,
$b^{-1}ab$ is in $H$.
There exists $u \in \left( {\Bbb Z}/s{\Bbb Z} \right)^\times$
such that $b^{-1}ab = a^u$.

Let $\{ \lambda_1, \lambda_2, \ldots, \lambda_n \}$ be the set of the eigenvalues of $a$,
where each $\lambda_i$ is a primitive $s$th root of $1$.
We regard it as a multi-set.

Then, by a famous theorem due to Frobenius,
$\{ \lambda_1^u, \lambda_2^u, \ldots, \lambda_n^u \}$ is the set of eigenvalues of $a^u$.

Since $b^{-1}ab = a^u$, 
\[
\{ \lambda_1, \lambda_2, \ldots, \lambda_n \}
= \{ \lambda_1^u, \lambda_2^u, \ldots, \lambda_n^u \}
\]
is satisfied as a multi-set.
Repeating it, we have
\begin{equation}\label{multi-set}
\{ \lambda_1, \lambda_2, \ldots, \lambda_n \}
= \{ \lambda_1^{(u^m)}, \lambda_2^{(u^m)}, \ldots, 
\lambda_n^{(u^m)} \}
\end{equation}
as a multi-set for any positive integer $m$.
Let ${\rm ord}(u)$ be the order of $u$ in the multiplicative group 
$\left( {\Bbb Z}/s{\Bbb Z} \right)^\times$.
Then, for any $i$,
\begin{equation}\label{eigen}
\left\{
\lambda_i, \lambda_i^u, \lambda_i^{(u^2)},  \ldots, \lambda_i^{(u^{{\rm ord}(u) - 1})}
\right\}
\end{equation}
is a subset of mutually distinct eigenvalues of the matrix $a$.
By (\ref{multi-set}), we know that eigenvalues in (\ref{eigen}) have 
the same multiplicity.
Therefore, it is easy to see that ${\rm ord}(u)$ divides $n$.
Since $n$ is a prime number, ${\rm ord}(u)$ is equal to $1$ or $n$.

\begin{itemize}
\item[(\rmn{1})]
If $u = 1$, then $G$ is abelian since $ab = ba$.
In this case, $G$ is cyclic as we have already seen in Subsection~\ref{abel}.
\item[(\rmn{2})]
Suppose ${\rm ord}(u) = n$.
Then, we may assume that 
\[
\{ \lambda, \lambda^u, \lambda^{(u^2)}, \ldots, \lambda^{(u^{n-1})} \}
\]
is the set of eigenvalues of $a$, where $\lambda$ is 
a primitive $s$th root of $1$.
Here, remark that the multiplicity of each eigenvalue is one.

Then there exists $c \in \gl{n}{K}$ such that
\begin{eqnarray}\label{matl}
c^{-1}ac=\left(
\begin{array}{ccccc}
\lambda &&&&O \\
&\lambda^u&&& \\
&&\lambda^{(u^2)}&& \\
&&&\ddots& \\
O&&&&\lambda^{(u^{n-1})}
\end{array}
\right) .
\end{eqnarray}
Replacing $G$ with $c^{-1}Gc$,
we may assume that $a$ is equal to the right-hand-side of (\ref{matl}).
Then,
\[
b^{-1}ab = a^u =
\left(
\begin{array}{ccccc}
\lambda^u &&&&O \\
&\lambda^{(u^2)}&&& \\
&&\ddots&& \\
&&&\lambda^{(u^{n-1})}& \\
O&&&&\lambda
\end{array}
\right) .
\]
By the above equality, the matrix $b$ coincides with
\[
({\bf b}_1 \ {\bf b}_2 \ \cdots \ {\bf b}_{n-1} \ {\bf b}_0) ,
\]
where ${\bf b}_i$ is an eigenvector of $a$ 
of eigenvalue $\lambda^{(u^i)}$ for $i = 0, 1, \ldots, n-1$.
Therefore, we may assume that the matrix $b$ is of the following form:
\[
\left( 
\begin{array}{ccccc}
0&\cdots&\cdots&0&b_0 \\
b_1&0&\cdots&\cdots&0 \\
0&b_2&\ddots&\ddots&\vdots \\
\vdots&\ddots&\ddots&\ddots&\vdots \\
0&\cdots&0&b_{n-1}&0
\end{array}
\right) 
\]
Then,
\[
\det(b) = (-1)^{n-1}b_0b_1\cdots b_{n-1} = 1 .
\]
On the other hand,
\begin{eqnarray*}
& & 
\det(te-b)  \\
& = &
\det \left( 
\begin{array}{ccccc}
t&0&\cdots&0&-b_0 \\
-b_1&t&\ddots&\ddots&0 \\
0&-b_2&t&\ddots&\vdots \\
\vdots&\ddots&\ddots&\ddots&0 \\
0&\cdots&0&-b_{n-1}&t
\end{array}
\right)  
\\
& = & \det \left( 
\begin{array}{ccccc}
	t&0&\cdots&\cdots&0 \\
	-b_1&t&0&\cdots&\vdots \\
	0&-b_2&t&\ddots&\vdots \\
	\vdots&\ddots&\ddots &\ddots &0 \\
	0&\cdots&0&-b_{n-1}&t
\end{array}
\right) + \det \left( 
\begin{array}{ccccc}
	0&0&\cdots&\cdots&-b_0 \\
	-b_1&t&0&\cdots&\vdots \\
	0&-b_2&t&\ddots&\vdots \\
	\vdots&\ddots&\ddots &\ddots &0 \\
	0&\cdots&0&-b_{n-1}&t
\end{array}
\right)
\\
& = & t^n+ (-1)^{n + (n-1)}b_0 b_1 \cdots b_{n-1} \\
& = & t^n+ (-1)^n .
\end{eqnarray*}
Since $n$ is an odd number,
we know that $1$ is an eigenvalue of the matrix $b$.
It is a contradiction.
Therefore, ${\rm ord}(u)$ is not $n$.
\end{itemize}
We have completed a proof in the case where $G$ is solvable.

\subsection{Final step in our proof of Lemma~\ref{main}}

In this subsection, we prove Lemma~\ref{main} 
without any additional assumptions.

Let $G$ be a group satisfying the assumption of Lemma~\ref{main}.
We prove Lemma~\ref{main} by induction on ${}^\# G$.
By induction, any proper subgroup of $G$ is cyclic.
Let $S_p$ be a $p$-Sylow subgroup of $G$ for each prime number $p$.

First, assume that $S_p$ is a normal subgroup of $G$ for any prime number $p$.
Then it is well known that $G$ is isomorphic to the direct product
of all Sylow subgroups.
Therefore, in this case, $G$ is solvable.
Thus, $G$ is cyclic as we have already seen in Subsection~\ref{skakai}.

Next, we assume that there exists a prime number $p$ such that
$S_p$ is not a normal subgroup of $G$.
The following subgroup is called  the {\em normalizer} of $S_p$.
\[
N_G(S_p) = \{
c \in G \mid cS_pc^{-1} = S_p \} 
\]
Since $S_p$ is not a normal subgroup of $G$,
$G \neq N_G(S_p)$.

Remember the following famous theorem due to 
Burnside (cf.\ Theorem~7.50 in \cite{Burnside}):

\begin{Theorem}[Burnside]
Let $F$ be a finite group.
Assume that there exists a prime number $q$ such that
a $q$-Sylow subgroup $S_q$ of $F$ is contained in 
the center of its normalizer $N_F(S_q)$.

Then there exists a normal subgroup $H$ of $F$
such that 
$$F=HS_q \ \ \mbox{and} \ \ H\cap S_q= \{ e \} .$$
\end{Theorem}

In our case, $S_p$ is contained in the center of $N_G(S_p)$
because $N_G(S_p)$ is cyclic.
By the above theorem due to Burnside,
there exists a normal subgroup $H$ of $G$
such that 
$$G=HS_p \ \ \mbox{and} \ \ H\cap S_p= \{ e \} .$$
Since $S_p \neq \{ e \}$, $H$ is a proper subgroup of $G$.
Therefore, $H$ is cyclic.
Since $S_p$ is a proper subgroup of $G$, $S_p$ is also cyclic.
Then, $G$ is solvable because of
\[
G/H \simeq S_p .
\]
Since $G$ is solvable, it is a cyclic group
as we have already seen in Subsection~\ref{skakai}.

We have completed a proof of Lemma~\ref{main}.

\section{The case where $n$ is not an odd prime number}\label{rei}
Suppose that $n$ is an integer bigger than $1$.

In this section, we give examples of non-abelian finite subgroups of 
$\sll{n}{{\Bbb C}}$ that satisfy the assumption in Lemma~\ref{main}
except for that $n$ is an odd prime number.

These examples are of type~\RMN{1} of  Theorem~6.1.11
and the representations are given in Theorem~5.5.6 in \cite{wolf}. 

\subsection{The case where $n$ is an even number}
In this subsection, we assume that $n$ is an even number.

Let $H$ be a non-abelian finite subgroup of $\sll{2}{{\Bbb C}}$.
For example, $H = \seisei{A, B}$, where
$$A=
\left(
	\begin{array}{cc}
		i &0\\
		0 &-i
	\end{array}
\right)
, \ \ B=
\left(
	\begin{array}{rr}
		0&i\\
		i&0
	\end{array}
\right) .
$$
It is easy to see that
$1$ is not an eigenvalue of any matrix in $H$ except for $e$.

Here we define as
\[
G = \left\{ \left.
\left(
\begin{array}{ccccc}
M      & 0 & \cdots & \cdots & 0 \\
0      & M &	0   & \cdots & 0 \\
\vdots & 0 & \ddots & \ddots & \vdots \\
\vdots & \vdots & \ddots & \ddots & 0 \\
0 & \cdots & \cdots & 0 & M
\end{array}
\right) \in \sll{n}{{\Bbb C}} \
\right| \
M \in H
\right\} .
\]
Then $1$ is not an eigenvalue of 
any matrix in $G$ except for $e$.
Since $G$ is isomorphic to $H$ as a group,
$G$ is not abelian.

\subsection{The case where $n$ is an odd composite number}
In this subsection,
assume that $n$ is an odd composite number.

Set 
\begin{equation}\label{n=qn'}
n = qn', 
\end{equation}
where $q$ is an odd prime number and $n'$ is an odd number
such that $q\leq n'$.

By a famous theorem due to Dirichlet, there exists an odd prime number $l$
such that
\begin{eqnarray*}
l\equiv 1 \pmod{2q} .
\end{eqnarray*}
Then, there exists $\alpha \in ({\Bbb Z}/l{\Bbb Z})^\times$
such that the order of $\alpha$ is $q$, i.e., it satisfies
\begin{eqnarray}\label{41alpha}
\alpha^q \equiv 1 \pmod{l} \ \ \ \mbox{and} \ \ \
\alpha \not\equiv \ \ 1 \pmod{l} .
\end{eqnarray}
Let $z$ (resp.\ $x$) be a primitive 
$l$th root (resp.\ $q$th root) of $1$.

Here, set
$$
A=
\left(
	\begin{array}{ccc|c}
		&O&&x\\\hline
		1&&O&\\
		&\ddots &&O\\
		O&&1&
	\end{array}
\right) 
, \ \
B=\left(
\begin{array}{cccc}
	z& & & O\\
	 &z^\alpha & & \\
	 & &\ddots & \\
	O& & &z^{(\alpha^{q-1})}
\end{array}
\right)  \in \gl{q}{{\Bbb C}}.
$$
\begin{Lemma}\label{4.3}
Set $G=\seisei{A,B} \subset \gl{q}{{\Bbb C}}$.
Then we have the following:
\begin{description}
\item[(\rmn{1})]
$\det A=x$, $\det B=1$.
\item[(\rmn{2})]
$AB\neq BA$. In particular, $G$ is not abelian.
\item[(\rmn{3})]
$G$ is a finite group.
\item[(\rmn{4})]
$1$ is not an eigenvalue of any matrix in $G$ except for the unit matrix.
\end{description}
\end{Lemma}

\proof
We have
\begin{eqnarray*}
\det A & = & 
(-1)^{q-1}x=x \\
\det B & = &
\prod_{i=0}^{q-1} z^{(\alpha^i)}=z^{\frac{\alpha^q-1}{\alpha -1}} .
\end{eqnarray*}
Since $l$ divides $\frac{\alpha^q-1}{\alpha -1}$ by (\ref{41alpha}),
\begin{eqnarray*}
z^{\frac{\alpha^q-1}{\alpha -1}}=1 .
\end{eqnarray*}
The assertion (\rmn{1}) has been proven.

\begin{eqnarray*}
A^{-1}BA & = &
\left(
	\begin{array}{c|ccc}
		&1&&O \\
		O&&\ddots & \\
		&O&&1 \\ \hline
		x^{-1}&&O&
	\end{array}
\right)
\left(
\begin{array}{cccc}
	z& & & O \\
	 &z^\alpha & & \\
	 & &\ddots & \\
	O& & &z^{(\alpha^{q-1})}
\end{array}
\right)
\left(
	\begin{array}{ccc|c}
		&O&&x \\ \hline
		1&&O& \\
		&\ddots &&O \\
		O&&1&
	\end{array}
\right) \nonumber \\
&=&
\left(
\begin{array}{ccccc}
	z^\alpha& & & & O \\
	 &z^{(\alpha^2)}& & & \\
	 & &\ddots & & \\
	 & & &z^{(\alpha^{q-1})} & \\
	O& & & &z
\end{array}
\right) =B^\alpha
\end{eqnarray*}
Since $z \neq z^\alpha$, we have $AB\neq BA$.
The assertion (\rmn{2}) has been proven.

It is easy to see that the order of $B$ is $l$.
Since
$$A^q=\left( \begin{array}{ccc}
	x& & O\\
	 &\ddots & \\
	O& &x
\end{array}
\right) ,$$
the order of $A$ is $q^2$.
Since $BA=AB^\alpha$, we have
\begin{eqnarray*}
G=\{ A^rB^s | r=0,1,\ldots ,q^2-1; \ s=0,1,\ldots ,l-1\} .
\end{eqnarray*}
In particular, the order of $G$ is finite.
The assertion (\rmn{3}) has been proven.

Now, we shall show that $1$ is not an eigenvalue of $A^rB^s$
for $r=0,1,\ldots ,q^2-1$, $s=0,1,\ldots ,l-1$
except for the case $r = s = 0$.

Set 
\begin{eqnarray*}
r=uq+v ,
\end{eqnarray*}
where $u$ and $v$ are integers such that $0 \leq u,  v < q$.

First, assume $v=0$.
Since
$$A^rB^s=x^u \left(
\begin{array}{cccc}
	z^s& & & O\\
	 &z^{s\alpha} & & \\
	 & &\ddots & \\
	O& & &z^{s\alpha^{q-1}}
\end{array}
\right)
=
\left(
\begin{array}{cccc}
	x^uz^s& & & O\\
	 &x^uz^{s\alpha} & & \\
	 & &\ddots & \\
	O& & &x^uz^{s\alpha^{q-1}}
\end{array}
\right) ,
$$
$\{ x^uz^s, \ x^uz^{s\alpha}, \ \ldots \ ,x^uz^{s\alpha^{q-1}} \}$ is
the set of the eigenvalues of $A^rB^s$.
Here assume that $x^uz^{s\alpha^t}=1$ for some $0 \leq t \leq q-1$.
Since $q$ and $l$ are relatively prime, we have
\begin{eqnarray*}
u\equiv 0 \pmod{q} \\
s\alpha^t\equiv 0 \pmod{l} .
\end{eqnarray*}
Therefore, we have $r=s=0$.

Next assume $v\neq 0$.
\begin{eqnarray*}
A^rB^s & = & (A^q)^uA^vB^s \\
&=& 
\begin{array}{l}
\hphantom{a}
\overbrace{
\hphantom{
\begin{array}{ccc}
x^u&\ddots&x^u
\end{array}
}}^{q-v}
\overbrace{
\hphantom{
\begin{array}{ccc}
x^{u+1}&\ddots&x^{u+1}
\end{array}
}}^{v}
\\
\left(
\begin{array}{ccc|ccc}
	&&&x^{u+1}& & 0 \\
	&O&& &\ddots & \\
	&&&0& &x^{u+1} \\\hline
	x^u&&0&&& \\
	&\ddots &&&O&\\
	0&&x^u&&&
\end{array}
\right)
\left(
\begin{array}{cccc}
	z^s& & & O\\
	 &z^{s\alpha} & & \\
	 & &\ddots & \\
	O& & &z^{s\alpha^{q-1}}
\end{array}
\right)
\end{array}  \\
&=& 
\left(
\begin{array}{ccc|ccc}
	&&&x^{u+1}z^{s\alpha^{q-v}}& & O \\
	&O&& &\ddots & \\
	&&&O& &x^{u+1}z^{s\alpha^{q-1}} \\\hline
	x^uz^s&&O&&& \\
	&\ddots &&&O&\\
	O&&x^uz^{s\alpha^{q-v-1}}&&&
\end{array}
\right) .
\end{eqnarray*}
Therefore, we know that
$$\mbox{the $(i,j)$th entry of $tE-A^rB^s$}=
\left\{ \begin{array}{ll}
t&(i=j) \\
-x^uz^{s\alpha^{j-1}} & (i=j+v) \\
-x^{u+1}z^{s\alpha^{j-1}} & (i=j+v-q) \\
0 & (\mbox{otherwise}) .
\end{array}
\right.$$
For each $j$, the $(i,j)$th entry of $tE - A^rB^s$ is not $0$
if and only if $i=j$ or $i \equiv j+v \pmod{q}$.
Since $q$ and $v$ are relatively prime,
we have
\begin{eqnarray*}
\det (tE-A^rB^s)&=&
t^q+(-1)^{q+v(q-v)} x^{uq+v}z^{s(1+\alpha + \cdots +\alpha^{q-1})}\\
&=&t^q-x^v .
\end{eqnarray*}
Since $x^v \neq 1$,
$1$ is not an eigenvalue of $A^rB^s$.
\qed

We define a group homomorphism
\[
f : G \longrightarrow \gl{qn'}{{\Bbb C}}
\]
by 
\[
\begin{array}{rccc}
&&&\overbrace{
\hphantom{
\begin{array}{ccc}
	C&\ddots&O
\end{array}
} }^{\frac{q+n'}{2}} 
\hphantom{
\begin{array}{ccc}
	C&\ddots&O
\end{array}
} \\
&f(C)& = &
\left(
\begin{array}{ccc|ccc}
	C&&O&&&\\
	&\ddots &&&O&\\
	O&&C&&&\\\hline
	&&&\overline{C}&&O\\
	&O&&&\ddots &\\
	&&&O&&\overline{C} 
\end{array}
\right)\\
&&&
\hphantom{
\begin{array}{ccc}
	C&\ddots&O
\end{array}
}
\underbrace{
\hphantom{
\begin{array}{ccc}
	C&\ddots&O
\end{array}
}}_{\frac{n'-q}{2}}
\end{array}
\]
for each $C \in G$, where
$\overline{C}$ is the complex conjugate matrix of $C$.
Here, remember that $n'$ is an odd number satisfying (\ref{n=qn'}).
If $C$ is not the unit matrix, 
$1$ is an eigenvalue of neither $C$ nor $\overline{C}$.
Therefore, if $C$ is not the unit matrix, 
$1$ is not an eigenvalue of $f(C)$.

On the other hand,
\begin{eqnarray*}
\det f(A)=(\det A)^{\frac{q+n'}{2}}(\det \bar{A})^{\frac{n'-q}{2}}
=x^{\frac{q+n'}{2}}(x^{-1})^{\frac{n'-q}{2}}
=x^q=1
\end{eqnarray*}
and, obviously $\det f(B)=1$.
Therefore, $f(G)\subset \sll{n}{\Bbb C}$.
Since $AB\neq BA$, 
\[
f(A)f(B)\neq f(B)f(A) .
\]
Therefore, $f(G)$ is not abelian.

\noindent
\begin{tabular}{l}
Kazuhiko Kurano \\
Department of Mathematics \\
Faculty of Science and Technology \\
Meiji University \\
Higashimita 1-1-1, Tama-ku \\
Kawasaki 214-8571, Japan \\
{\tt kurano@isc.meiji.ac.jp} \\
{\tt http://www.math.meiji.ac.jp/\~{}kurano}
\end{tabular}

\vspace{2mm}

\noindent
\begin{tabular}{l}
Shougo Nishi \\
Department of Mathematics \\
Faculty of Science and Technology \\
Meiji University \\
Higashimita 1-1-1, Tama-ku \\
Kawasaki 214-8571, Japan \\

\end{tabular}


\begin{thebibliography}{99}

\bibitem{Benson}
{\sc D. J. Benson},
{\em Polynomial invariants of finite groups},
London Mathematical Society Lecture Note Series {\bf 190},
Cambridge University Press, Cambridge, 1993.

\bibitem{Nagata}
{\sc M. Nagata}, 
{\em Local rings},
Interscience Tracts in Pure and Applied Mathematics No.\ 13.
Interscience Publishers a division of John Wiley \& Sons, New York-London 1962.

\bibitem{Burnside}
{\sc J.\ J.\ Rotman}
{\em An introduction to the theory of groups, Fourth edition},
Graduate Texts in Mathematics {\bf 148}, Springer-Verlag, New York, 1995.

\bibitem{Watanabe}
{\sc K.-i. Watanabe},
{\em Certain invariant subrings are Gorenstein \RMN{1}, \RMN{2}},
Osaka J. Math.\ {\bf 11} (1974), 1--8, 379--388.

\bibitem{Watanabe1}
{\sc K.-i. Watanabe},
{\em Invariant theory of finite groups (in Japanese)},
Progress in group theory, pp135--183,
Asakura 2004.

\bibitem{wolf}
{\sc J. A. Wolf},
{\em Spaces of constant curvature},
fifth edition, Publish or Perish, Inc.\ 1984.

\bibitem{YY}
{\sc S. S.-T. Yau} and {\sc Y. Yu},
{\em Gorenstein quotient singularities in dimension three},
Mem. Amer. Math. Soc. Vol.\ 105 (1993), No\ 505.
\end{thebibliography}
\end{document}